\documentclass[journal]{IEEEtran}
\usepackage{cite}
\usepackage{amsmath,amssymb,amsfonts}
\usepackage{graphicx}
\usepackage{textcomp}
\usepackage{float}
\usepackage[colorlinks=true,linkcolor=magenta,citecolor=blue,urlcolor=blue,]{hyperref}
\usepackage{balance}
\usepackage{bm,bbm,mathrsfs,amssymb,pifont,amssymb,pifont,amsfonts,amsmath,stmaryrd,color,amssymb,graphics,graphicx,epstopdf}

\newtheorem{theorem}{\bf Theorem}
\newtheorem{assumption}{\bf Assumption}
\newtheorem{problem}{\bf Problem}
\newtheorem{remark}{\bf Remark}

\newtheorem{lemma}{\bf Lemma}

\newcommand{\sgn}{\mbox{sgn}}
\newcommand*{\QEDB}{\hfill\ensuremath{\square}}

\ifCLASSINFOpdf
 \else
 \fi

\begin{document}

\title{Nonlinear Bipartite Output Regulation with Application to {T}uring Pattern}

\author{Dong~Liang, Martin~Guay, Shimin~Wang\thanks{This work has been partly supported by the National Natural Science Foundation of China under Grant  62203305 and Natural Sciences and Engineering Research Council (NSERC), Canada. (Corresponding author: Shimin Wang.)
	}
\thanks{D. Liang is with the Department of Control Science and Engineering, School of Optical-Electrical and Computer Engineering, University of Shanghai for Science and Technology, Shanghai 200093, China. (E-mail: dliang@usst.edu.cn);
Martin~Guay and Shimin~Wang are with the Department of Chemical Engineering, Queen's University, Kingston,  ON K7L 3N6, Canada (E-mail: martin.guay@queensu.ca, shimin.wang@queensu.ca)}}

\maketitle

\begin{abstract}
In this paper, a bipartite output regulation problem is solved for a class of nonlinear multi-agent systems subject to static signed communication networks. A nonlinear distributed observer is proposed for a nonlinear exosystem with cooperation-competition interactions to address the problem. Sufficient conditions are provided to guarantee its existence and stability. The exponential stability of the observer is established. As a practical application, a leader-following bipartite consensus problem is solved for a class of nonlinear multi-agent systems based on the observer. Finally, a network of multiple pendulum systems is treated to support the feasibility of the proposed design. The possible application of the approach to generate specific {T}uring patterns is also presented.
\end{abstract}

\begin{IEEEkeywords}
Distributed observer, output regulation, bipartite consensus, nonlinear multi-agent systems, signed graphs.
\end{IEEEkeywords}

\IEEEpeerreviewmaketitle

\section{Introduction}

\IEEEPARstart{O}ver the past few decades, the design of distributed control law of multi-agent systems (MASs) subject to cooperation-competition interactions modelled by the signed graphs has gathered significant attention within the system and control community.
A signed communication graph can be treated as an extension of the traditional unsigned communication graphs where the weight of the signed situation can be negative, representing the competitive relationships.
A typical problem related to this field is the so-called bipartite consensus problem. This problem aims to design a distributed controller for MASs with cooperation-competition interactions so that the bipartite property can be achieved \cite{altafini2013consensus}.
Under the condition that the signed communication graph is structurally balanced, the nodes of the MASs can be partitioned into two sets based on the sign of the edge weights.
In general, cooperation and competition relationships cannot live independently. Signed communication graphs are, therefore very relevant in real-world applications. They can be used to model the relationships of several scenarios, such as allied/adversary relationships, edge dynamics \cite{du2019edge,du2019edge2}, systems biology \cite{griffin2004cooperation}, and the design of the opinion dynamics in social networks \cite{altafini2015predictable,Meo15,Xia16}, to name a few.

The bipartite consensus problem was first investigated in \cite{altafini2013consensus} as an extension of the classical consensus problem.
This result, limited to the bipartite consensus problem for single-integrator MASs, has been generalized to more general linear MASs \cite{yaghmaie2017bipartite,Zha17,zhang2019h}.
The so-called {\em gauge transformation} was defined in \cite{altafini2013consensus} to convert the closed-loop system from the signed communication graph to the unsigned case.
The bipartite consensus problem can be sorted into the leaderless case
\cite{altafini2013consensus} and the leader-following case  \cite{Jia18,Ma18,bhowmick2018leader}.
For a MAS consisting of a leader node, different bipartite tracking problems have been investigated in \cite{shi2023cooperation,Jia18, Ma18, Liu18,shao2020}.
For example, a leader system with scalar dynamics was proposed in \cite{Ma18}, while a leader system with a high dimension has been adopted in \cite{Jia18}.
Robustness and matching uncertainty have been considered in \cite{Liu18} while authors in \cite{Ma18} considered the measurement noise.
Their analysis provided sufficient and necessary conditions to achieve the bipartite tracking property. In addition, the related sign consensus problem has been investigated in \cite{jiang2017sign,jiang2020output,sun2021output}.
Static signed communication graphs were considered in \cite{jiang2017sign} while switching topologies were addressed in \cite{jiang2020output}.

More recently, some researchers have focused on the design of output-based distributed control law to settle the leader-following bipartite synchronization problem of linear MASs \cite{Wan17}.
In particular, the bipartite output regulation problem (BORP) has been considered as early as in \cite{aghbolagh2017bipartite}.
The BORP can be treated as an extension of the traditional cooperative output regulation problem of linear MASs subject to unsigned graphs studied in \cite{shtac12} to the signed case.
A new distributed observer over static signed communication graphs was proposed in \cite{aghbolagh2017bipartite} for systems subject to a linear exosystem.
In \cite{liang2020robust}, a robust BORP approach was proposed for general linear uncertain MAs in which state-based and output-based distributed control laws were proposed.
It should be noted that current results on the BORP remain limited to linear dynamical systems. In \cite{shams21hinfinity}, the bipartite tracking consensus problem over a directed signed graph was studied for a class of nonlinear MASs in a particular form where the nonlinear dynamics and external disturbances are added to the linear dynamics of each follower. Other recent advances in nonlinear MASs can be found in \cite{Qin17} in which the input saturation was added to the general linear MASs and in \cite{liang2020leader} for Euler-Lagrange systems.  In \cite{Zha16}, a type of nonlinear MASs was investigated, achieving bipartite synchronization subject to switching topologies. Nevertheless, the nonlinear systems considered in \cite{Zha16} were 1-dimension satisfying the one-sided Lipschitz condition and the problem in \cite{Zha16} is leaderless without considering the leader node. In this technical paper, we further investigate a type of leader-follower bipartite synchronization problem denoted by the nonlinear BORP. From the literature, most of the research on bipartite synchronization is restricted to linear MASs. However, designing a distributed control law for the bipartite consensus problem for general nonlinear MASs is still challenging. The reason is that the {\em gauge transformation} defined in \cite{altafini2013consensus} only applies to linear systems, and it loses efficacy if the closed-loop system is nonlinear. Developing ``nonlinear" gauge transformation techniques for nonlinear MASs remains an open problem.

The pioneering work on the {T}uring pattern dates back to the result in \cite{turing1990}, where {T}uring proposed the fundamental mathematical theory of formulating the pattern based on morphogenesis. Actually, various Turing patterns exist widely in nature and society, e.g., system biology, neuroscience, chemistry, physics, ecology, etc. \cite{MaronSusan1997,maini2006turing}. Typical {T}uring patterns can be found in zebra stripes, the skin of fish, spirals in chemical processes, and so on \cite{nakao2010turing}. More recently, dynamical systems such as reaction–diffusion systems have been employed to investigate the formulation of {T}uring patterns subject to undirected or directed communication graphs \cite{asllani2014theory}. In this paper, we try to relate signed graphs, which can be treated as the generalization of directed graphs, with the {T}uring pattern arising from the networked nonlinear dynamical systems. 
Furthermore, we formulate the nonlinear BORP of nonlinear MASs over signed communication graphs. The nonlinear BORP can be treated as a generalization of the traditional cooperative output regulation problem (CORP) as studied in \cite{shtac12,Su14,liu2019distributed}. The research on the CORP for general linear MASs appeared in \cite{shtac12,Su14}. The result was generalized to systems subject to switching topologies in \cite{su2012cooperative}.

It should be noted that \cite{huang2017cooperative} and \cite{wang2018adaptive} proposed two different kinds of adaptive distributed observers to deliver and estimate the leader's known and unknown system matrix, respectively.
Moreover, the distributed observer approach has been found in many applications to control classes of nonlinear MASs. Examples include Euler-Lagrange systems \cite{wang2018adaptive},  rigid spacecraft systems \cite{wang2021event}, amongst others. Other recent advances can be found in \cite{liu2019distributed}.
To the best of our knowledge, the nonlinear BORP has not been considered in the literature. In nonlinear BORP, asymptotic tracking and disturbance rejection are considered for a class of nonlinear MASs where a nonlinear autonomous system generates the tracking reference signal called the exosystem.
This nonlinear exosystem can produce a large class of standard signals such as steps, ramps,  sinusoidal, specified nonlinear functions or their combination. To deal with the BORP, we propose a nonlinear distributed observer (NDO) over static signed communication graphs.
A technical lemma is then established, and sufficient conditions are given to ensure the existence and the exponential stability of the NDO. We show that, under the assumption that the signed communication digraph is structurally balanced and contains a spanning tree with the leader node as the root, the NDO can estimate the original and the opposite state of the nonlinear leader system exponentially.
In this design, the state of some agents synchronizes with the leader's state, while others are forced to follow the opposite of the leader's state. As opposed to the results in \cite{liu2019distributed,Su14,su2012cooperative}, which focused on the unsigned graphs situation, the proposed approach applies to more general static signed graphs.
The authors in \cite{aghbolagh2017bipartite} extended the distributed observer for linear exosystem from the unsigned communication graphs to the signed case. In this technical paper, a nonlinear exosystem is investigated. Since nonlinear terms exist in the NDO, the {\em gauge transformation} defined in \cite{altafini2013consensus} cannot be directly implemented for the nonlinear systems.
To alleviate this difficulty, we employ a class of {\em ``nonlinear gauge transformation''} to convert the nonlinear closed-loop system from the signed case to the unsigned case, which bridges the signed and unsigned cases.
As an application, a state feedback dynamic distributed control law is further synthesized to solve the leader-following bipartite consensus problem for a class of nonlinear MASs based on the NDO via the certainty equivalence principle.

The main contribution of this study can be summarized as follows:

1) A NDO is proposed for a class of nonlinear autonomous exosystems for systems operated over static signed graphs. The NDO in \cite{liu2019distributed} is extended from the unsigned communication graphs to the signed situation. Compared to existing results on the distributed observer \cite{shtac12,aghbolagh2017bipartite,liu2019distributed}, the proposed design is more general and applies to nonlinear autonomous exosystems with cooperation-competition interactions.

2) A technical lemma is established, and sufficient conditions are provided to guarantee the existence of the NDO over static signed communication graphs. The concept of ``nonlinear gauge transformation" is introduced to convert the nonlinear closed-loop system from a signed case to an unsigned case which bridges the design of distributed control law over signed and unsigned scenarios.

3) On the basis of the NDO and the certainty equivalence principle, a distributed state feedback control law is designed to solve the leader-following bipartite consensus problem for a class of nonlinear MASs. The design technique implemented in this paper can also be used to design distributed control law to solve the bipartite synchronization problem of other nonlinear MASs.  Moreover, we apply our NDO to drive multiple networked nonlinear dynamical systems to behave in a Turing pattern similar to the zebra stripe.

The remaining part of this paper is organized in the following order. Section \ref{pre} briefly summarizes the basic definition of the signed communication graphs and related lemmas. The nonlinear BORP is formulated in this section formally.  Our main result on the NDO will be illustrated in Section \ref{ndo}. Section \ref{application} formulates and solves the leader-following bipartite consensus for a class of nonlinear MASs. A simulation study is presented in Section \ref{exam} where the application to the {T}uring pattern is addressed. Some conclusions follow this in Section \ref{con}.

{\bf Notation:} For a $n$-dimensional vector $z \in \mathbb{R}^n$, $||z||$ denotes the 2-norm. The operator $\otimes$ is the Kronecker product of matrices. The symbol $\mathbf{1}_N$ represents an $N$-dimensional column vector whose elements are $1$. For a column vector $a_i$, $i=1,2,\dots,s$, $\mbox{col} (a_1,\dots,a_s)=[a_1^T,\dots,a_s^T]^T$. $\textnormal{diag}(x_1,x_2,\dots,x_n)$ represents a diagonal matrix whose diagonal element is $x_1$, $x_2$, $\dots$, $x_n$.  $\sgn(x)$ denotes the signum function of a real number $x \in \mathbb{R}$ and $\sgn(x)=-1$ if $x<0$, $\sgn(x)=0$ if $x=0$ and $\sgn(x)=1$ if $x>0$.

\section{Preliminaries and Problem Formulation}\label{pre}
\subsection{Signed Graphs}
A signed communication graph is denoted by $\mathcal{G}^s = (\mathcal{V}, \mathcal{E}, \mathcal{A}^s)$ where $\mathcal{V}=\{1,2,\dots,N\}$ represents the set of nodes, $\mathcal{E} \subseteq \mathcal{V} \times \mathcal{V}$ is the set of edges, and $\mathcal{A}^{s}=[a_{ij}] \in \mathbb{R}^{N\times N}$ is the adjacency matrix of the signed communication graph $\mathcal{G}^s$.
The notation $(j, i)$ with $j \neq i$ denotes an edge of $\mathcal{E}$ from the node $j$ to the node $i$.
 The elements $a_{ij}$, $i, j = 1, \cdots, N$,  of the matrix $\mathcal{A}^{s}$ are such that  $a_{ii}=0$, and,  for $i \neq j$,   $a_{ij} \neq 0$  iff $(j, i) \in \mathcal{E}$. Let $\mathcal{N}_i = \{j \; |\; (j, i) \in \mathcal{E}\}$, which is defined by the neighboring set of the node $i$.
 If there exists an edge set $(i_1, i_2), (i_2, i_3), \dots, (i_{k-1}, i_{k})$ in a directed graph (digraph), this edge set is called a directed path from the node $i_1$ to the node $i_k$. In this kind of situation, node $i_{k}$ is said to be reachable from the initial node $i_1$. A digraph is said to contain a spanning tree if there is a node $i$ reaching all other nodes. In this case,  the node $i$ is said to be the root node. The Laplacian matrix $\mathcal{L}$ of a signed communication graph $\mathcal{G}^s$ is defined by $\mathcal{L}^s=\mathcal{C}-\mathcal{A}^s$ where the degree matrix $\mathcal{C}=\textnormal{diag}(c_1,\dots,c_N)$ is diagonal and $[\mathcal{C}]_i = \sum_{j \in \mathcal{N}_i} |a_{ij}|$.  A signed communication graph is structurally balanced if it can be partitioned into two disjoint sets of nodes $\mathcal{V}_1$ and $\mathcal{V}_2$ with $\mathcal{V}_1\cup\mathcal{V}_2=\mathcal{V}$ and $\mathcal{V}_1\cap\mathcal{V}_2=\varnothing$, such that, $a_{ij} \geq 0$, $\forall i, j \in \mathcal{V}_q$,  where $q \in \{1,2\}$,  $a_{ij}\leq 0$, $\forall i \in \mathcal{V}_p$, $\forall j \in \mathcal{V}_q$, $p \neq q$, where $p, q \in \{1,2\}$.
It is called structurally unbalanced otherwise.

\subsection{Problem Formulation}
In this part, we formally formulate the nonlinear BORP of nonlinear MASs over static signed communication graphs.

Let us consider a type of general nonlinear MASs described as follows:
\begin{eqnarray} \label{plant0}
\begin{split}
\dot{x}_i&=f_i(x_i,u_i,v)\\
y_i&=h_i(x_i,u_i,v)\\
e_i&=y_i -\phi_iy_0, \,\, i= 1,\dots,N
\end{split}
\end{eqnarray}
where $x_i \in \mathbb{R}^{n_i}$ is the state for agent $i$, $u_i \in \mathbb{R}^{p_i}$ denotes the control input for agent $i$, $y_i \in \mathbb{R}^{q}$ represents the measurement output and  $e_i \in \mathbb{R}^{q}$ denotes the tracking error of each agent. $f_i(\cdot)$ and $h_i(\cdot)$ are some sufficiently smooth functions defined globally. $\phi_i \in \{1,-1\}$, $i=0,1,\dots,N$ and $\phi_i = 1$ if $i \in \mathcal{V}_1$, $\phi_i = -1$ if $i \in \mathcal{V}_2$. $v \in \mathbb{R}^{m}$ is the exogenous signals representing the tracking reference and/or external disturbances generated by the nonlinear exosystem:
\begin{eqnarray} \label{exo}
\begin{split}
\dot{v}&=a(v)\\
y_0&=g(v)
\end{split}
\end{eqnarray}
where $y_0 \in \mathbb{R}^{q}$ is the tracking reference signals. The functions $a(v)$ and $g(v)$ are defined globally and $a(0)=0$, $g(0)=0$.

The system consisting of the plant (\ref{plant0}) and the exosystem (\ref{exo}) can form a MAS with $N+1$ agents where the leader node $0$ is the exosystem (\ref{exo}).  For a signed digraph $\bar{\mathcal{G}}^s=(\bar{\mathcal{V}}, \bar{\mathcal{E}}, \bar{\mathcal{A}}^s)$ with a leader node $0$ (\ref{exo})  where $\bar{\mathcal{V}}=\{0,1,2,\dots,N\}$, the Laplacian matrix $\bar{\mathcal{L}}^s$ can be defined as follows:
\begin{eqnarray*}
	\bar{\mathcal{L}}^s=\left[ \begin{array}{c|c}
		\sum_{j=1}^N |a_{0j}|  & [a_{01},\dots,a_{0N}]\\ \hline
		-\Delta\mathbf{1}_N & \mathcal{H}^{s}\\
	\end{array} \right]
\end{eqnarray*}
where $\mathcal{H}^s=\mathcal{L}^s+\Delta$, $\mathcal{L}^s$ is the Laplacian matrix of
the signed communication digraph by removing the leader node and $\Delta=\textnormal{diag}(a_{10},\dots,a_{N0})$. It is assumed that $a_{i0} \geq 0,i=0,\dots,N$, w.l.o.g., meaning that the edge weight from the leader node to the follower is nonnegative.

We consider a class of distributed control laws of the form
\begin{eqnarray} \label{u}
\begin{split}
u_i&=k_i(x_i, \eta_i)\\
\dot{\eta}_i&= \theta_i(\eta_i, \eta_j-\eta_i, \phi_i)
\end{split}
\end{eqnarray}
where $k_i(\cdot)$ and $\theta_i(\cdot)$ are sufficiently smooth functions.

\begin{problem}[Bipartite Output Regulation of Nonlinear MASs]
	Given the nonlinear MASs (\ref{plant0}) with the nonlinear exosystem (\ref{exo}) and a signed communication digraph $\bar{\mathcal{G}}^s$, design a distributed control law $u_i$ of the form (\ref{u}) for each agent, such that\\
1) the trajectories of all states of the closed-loop system exist and are bound for all $t \geq 0$;\\
2) the tracking error for each agent satisfies
	\begin{eqnarray}
	\begin{split}
	\lim_{t \rightarrow \infty} e_i(t)=0,\,\,i=1,\dots,N.
	\end{split}
	\end{eqnarray}
\end{problem}

\begin{remark}
If the communication graph is unsigned, the nonlinear BORP of MASs defined in {\em Problem 1} reduces to the traditional nonlinear CORP defined in \cite{liu2019distributed}. The nonlinear BORP of MASs can also be viewed as  {\em the leader-following bipartite synchronization problem} for a class of nonlinear MASs where asymptotic tracking and disturbance rejection are considered concurrently.
\end{remark}

\begin{remark}
If property 2) in {\em Problem 1} is achieved, the leader-following bipartite synchronization is also achieved since the following property is satisfied:
\begin{eqnarray}
\begin{split}
&\lim_{t \rightarrow \infty} ||y_i(t) - y_0(t)||=0,\quad i \in \mathcal{V}_1\\
&\lim_{t \rightarrow \infty} ||y_i(t) + y_0(t)||=0, \quad i \in \mathcal{V}_2.
\end{split}
\end{eqnarray}
\end{remark}

\begin{remark}
To deal with {\em Problem 1},  a distributed control law can be synthesized on the basis of the distributed observer and purely decentralized control law, also called the certainty equivalence principle following \cite{shtac12,liu2019distributed}. Hence, the design procedure can be divided into two steps:\\
1) First, one designs an NDO over signed graphs for estimating the original and the opposite state for the nonlinear exosystem (\ref{exo});\\
2) Second, one designs a distributed control law (\ref{u}) on the basis of the NDO achieving the bipartite output regulation property.
\end{remark}

\subsection{Some Lemmas}
\begin{lemma}{(\cite{aghbolagh2017bipartite})}\label{l1}
For a signed communication digraph satisfying the structurally balanced property, there exists a diagonal matrix defined by the gauge transformation
\begin{eqnarray*}
\bar{\Phi}=\left[ \begin{array}{c | c}
   \phi_0  & 0\\ \hline
   0& \Phi\\
\end{array} \right]
\end{eqnarray*}
where $\Phi=\textnormal{diag}(\phi_i)$ for $i=1,\dots,N$ with $\phi_i \in \{1,-1\}$ for $i=0,1,\dots,N$ such that $\bar{\Phi}\bar{\mathcal{L}}^{s}\bar{\Phi}\mathbf{1}_N=\bar{\mathcal{L}}\mathbf{1}_N=0$. Define $\mathcal{H}=\Phi\mathcal{H}^{s}\Phi$ where the matrix $\mathcal{H}$ is defined in accord with the unsigned communication digraph case, we attain the following equality:
$\mathcal{H}\mathbf{1}_N = \Delta\mathbf{1}_N$.
\end{lemma}

\begin{lemma}(Theorem 2.5.2 in \cite{michel1977qualitative})
There exists a diagonal matrix $P=\textnormal{diag}(p_1,\dots,p_N)$ whose diagonal elements $p_i$, $i=1,\ldots,p_N$, are positive constants such that
\begin{eqnarray}
Q=P\mathcal{H}+\mathcal{H}^TP\nonumber
\end{eqnarray}
which is a symmetric and positive definite matrix where $\mathcal{H}$ is the matrix consisting of the last $N$ rows and the last $N$ columns of the Laplacian matrix $\bar{\mathcal{L}}$ for the unsigned communication graph.
\end{lemma}

\begin{lemma}(Lemma 11.1 in \cite{chen2015stabilization})\label{l3}
For any function $f(x,d):\mathbb{R}^{n}\times \mathbb{R}^{l} \rightarrow \mathbb{R}$ that is continuously differentiable and vanishing at origin satisfying $f(0,d)=0$, there exists a nonnegative smooth function $m(x,d)\geq 0$, such that
\begin{eqnarray}
|f(x,d)| \leq ||x||m(x,d).\nonumber
\end{eqnarray}
\end{lemma}

\section{Nonlinear Observer over Signed Graphs} \label{ndo}
In this section, we propose an NDO over static signed graphs and establish a technical lemma to guarantee the existence and exponential stability.

We start the development by listing some assumptions.
\begin{assumption} \label{A1}
For any initial states, the solution to (\ref{exo}) exists and is bounded for all $t \geq 0$.
\end{assumption}

\begin{assumption} \label{A2}
The nonlinear term $a(v)$ can be written as $a(v)=Mv+T(v)v$ where $M$ is the Jacobian of $a(v)$ at $v=0$ and $T(v)=\textnormal{diag}(d_1(v),\dots,d_m(v))$, $d_i(v)\leq 0$, $i=1,2,\dots,m$.
\end{assumption}

\begin{assumption} \label{A3}
The signed communication digraph $\bar{\mathcal{G}}^s$ is structurally balanced and contains a spanning tree where the root is the leader node 0.
\end{assumption}


\begin{remark}
Assumption \ref{A1} guarantees that, for any $v(0) \in \mathbb{V}_0$ where  $\mathbb{V}_0$ is a compact set, there exists a compact set $\mathbb{V}$ such that $v(t) \in \mathbb{V}$ for $t\in[0,\infty)$. Assumption \ref{A2} is made in \cite{liu2019distributed} to ensure the existence and exponential stability of the closed-loop for the NDO. Both Assumptions \ref{A1} and \ref{A2} are needed for the solvability of the nonlinear BORP. Assumption \ref{A3} is standard in the literature on  bipartite consensus problem \cite{altafini2013consensus} and MASs.
\end{remark}

The NDO over static signed graph is designed as follows:
\begin{eqnarray} \label{do2}
\begin{split}
\dot{\eta}_i=&\phi_ia(\phi_i\eta_i)+\mu\big(\sum_{j=1}^{N}a_{ij}(\eta_j-\sgn(a_{ij})\eta_i)\\
&\hspace{1in}+a_{i0}(\phi_iv-\eta_i)\big),
\end{split}
\end{eqnarray}
where $\eta_0=v$, $\eta_i \in \mathbb{R}^m$ and $\mu>0$ is some positive constant to be assigned. $\phi_i \in \{1,-1\}$ for $i=1,\dots,N$ and $\phi_i = 1$ if $i \in \mathcal{V}_1$, $\phi_i = -1$ if $i \in \mathcal{V}_2$.

\begin{remark}
The gauge transformation defined in \cite{altafini2013consensus} only applies to linear systems. The traditional gauge transformation cannot be employed directly in the nonlinear case. To circumvent this difficulty, we adopt expressions of the type $\phi_ia(\phi_i\eta_i)$ for the nonlinear function $a(\cdot)$, which can be treated as a nonlinear gauge transformation of the original closed-loop system to the unsigned case. This type of nonlinear gauge transformation will also be utilized in the design of distributed control law to solve the leader-following bipartite consensus problem of nonlinear MASs in Section \ref{application}.
\end{remark}

\begin{remark}
If the communication graph is unsigned, the NDO (\ref{do2}) reduces to the result in \cite{liu2019distributed}:
\begin{eqnarray}
\dot{\eta}_i=a(\eta_i)+\mu\sum_{j=0}^{N}a_{ij}(\eta_j-\eta_i).\nonumber
\end{eqnarray}
\end{remark}

\begin{remark}
If $a(v)$ is an odd function, i.e. $a(-v)=-a(v)$, (\ref{do2}) reduces to the following form:
\begin{eqnarray} \label{do}
\begin{split}
\dot{\eta}_i=&a(\eta_i)+\mu\big(\sum_{j=1}^{N}a_{ij}(\eta_j-\sgn(a_{ij})\eta_i)+a_{i0}(\phi_iv-\eta_i)\big).
\end{split}
\end{eqnarray}

Hence, the odd property is not necessary for the design of the NDO (\ref{do2}).
\end{remark}

\begin{remark}
If $a(v)=Sv$ where $S \in \mathbb{R}^{q \times q}$ is a constant square matrix, (\ref{do2})
reduces to the result in \cite{aghbolagh2017bipartite}:
\begin{eqnarray}
\begin{split}
\dot{\eta}_i&=S\eta_i+\mu\big(\sum_{j=1}^{N}a_{ij}(\eta_j-\sgn(a_{ij})\eta_i)+a_{i0}(\phi_iv-\eta_i)\big).
\end{split}\nonumber
\end{eqnarray}
\end{remark}

We then establish the following lemma for NDO (\ref{do2}).
\begin{lemma}\label{l5}
Consider the nonlinear exosystem (\ref{exo}), under Assumptions 1-3, for any compact subset $\mathbb{V}_0 \subset \mathbb{R}^m$ containing the origin, any initial conditions $v(0)$, $\eta_i(0)\in \mathbb{V}_0$, $i=1,2,\dots,N$, the solution to (\ref{do2}) exists for $t \geq 0$ and satisfies
\begin{eqnarray}
\lim_{t \rightarrow \infty} (\eta(t)-(\Phi\mathbf{1}_N) \otimes v(t)) = 0\nonumber
\end{eqnarray}
exponentially for $\mu>\mu_0$ where $\mu_0>0$ is a positive constant.
\end{lemma}

\emph{Proof:} Let $\eta=\textnormal{col}(\eta_1,\dots,\eta_N)$ and ${\bf a}(\eta)=\textnormal{col}(a(\eta_1),\dots,a(\eta_N))$.
Therefore, the NDO (\ref{do2}) can be rewritten into the following compact form:
\begin{eqnarray} \label{eta2}
\begin{split}
\dot{\eta}=&(\Phi \otimes I_m){\bf a}((\Phi \otimes I_m)\eta)-\mu(\mathcal{H}^s\otimes I_m)\eta \\
&+\mu(\Phi\Delta \otimes I_m)(\mathbf{1}_N \otimes v).
\end{split}
\end{eqnarray}


Let $z_i =\phi_i\eta_i$, $z=\textnormal{col}(z_1,\dots,z_N)$, $\Phi =\textnormal{diag}(\phi_1,\dots,\phi_N)$, ${\bf a}(z)=\textnormal{col}(a(z_1),\dots,a(z_N))$, the relationship between $z$ and $\eta$ is given as follows:
$z =(\Phi \otimes I_m)\eta$.
Moreover,
\begin{eqnarray} \label{closedz}
\begin{split}
\dot{z}
&=(\Phi \otimes I_m)((\Phi \otimes I_m){\bf a}(z)-\mu(\mathcal{H}^{s}\otimes I_m)\eta\\
&+\mu(\Phi\Delta \otimes I_m)(\mathbf{1}_N \otimes v))\\
&={\bf a}(z)-\mu(\mathcal{H}\otimes I_m)(z-\mathbf{1}_N \otimes v).
\end{split}
\end{eqnarray}

Hence, under Assumption \ref{A2}, the closed-loop system is given as follows:
\begin{eqnarray} \label{closed}
\begin{split}
\dot{z}=(I_N \otimes M+\mathcal{T}(z))z-\mu(\mathcal{H}\otimes I_m)(z-\mathbf{1}_N \otimes v)
\end{split}
\end{eqnarray}
where $\mathcal{T}(z)=\textnormal{diag}(T(z_1),\dots,T(z_N))$. Assume $v\equiv0$, then (\ref{closed}) reduces to
\begin{eqnarray} \label{z1}
\begin{split}
\dot{z}=(I_N \otimes M+\mathcal{T}(z))z-\mu(\mathcal{H}\otimes I_m)z.
\end{split}
\end{eqnarray}

Define the Lyapunov candidate function:
$V(z)=z^T(P \otimes I_m)z$,
which is positive definite. Under A.\ref{A3}, the eigenvalues of $\mathcal{H}$ have positive real parts and we attain $Q=P\mathcal{H}+\mathcal{H}^TP$ by Lemma 2.

The derivative along the trajectory of (\ref{z1}) is
\begin{eqnarray}
\begin{split}
\dot{V}(z) 
&=-\mu z^T(Q \otimes  I_m)z + 2z^T(P \otimes M)z\\
&+2z^T(P \otimes I_m)\mathcal{T}(z)z\\
&\leq -\mu \lambda_{min}(Q) ||z||^2+2\lambda_{max}(P \otimes M)||z||^2.
\end{split}\nonumber
\end{eqnarray}
Let $\mu_1=\frac{2\lambda_{max}(P \otimes M)}{ \lambda_{min}(Q)}$.
When $\mu > \mu_1$, $\dot{V}(z)$ is negative definite. Hence, we conclude that the origin of the system (\ref{z1}) is globally exponentially stable.

We then show that there exists a positive constant $\mu_2>0$ such that the origin of (\ref{closed}) is also exponential stable when $\mu>\mu_2$.

Therefore, there exists some positive constants $b_1$, $b_2$, $b_3$ and $b_4$ such that
\begin{eqnarray}
\begin{split}
&b_1||z||^2 \leq V(z) \leq b_2||z||^2\\
&\frac{\partial V(z)}{\partial z}({\bf a}(z)-\mu(\mathcal{H}\otimes I_m)z) \leq -b_3||z||^2\\
&\left|\left|\frac{\partial V(z)}{\partial z}\right|\right| \leq b_4 ||z||.
\end{split}\nonumber
\end{eqnarray}

The derivative of $V(z)$ along the trajectories of (\ref{closedz}) is given as follows:
\begin{eqnarray}
\begin{split}
\dot{V}(z) 
&=\frac{\partial V(z)}{\partial z}\Big({\bf a}(z)-\mu(\mathcal{H}\otimes I_m)z\Big)+\mu\frac{\partial V(z)}{\partial z}(\mathcal{H}\mathbf{1}_N\otimes v)\\
&\leq -b_3||z||^2+\mu b_4 ||z||||\mathcal{H}\mathbf{1}_N\otimes v||\\
& \leq -\frac{b_3-\lambda}{b_2}V(z)+\frac{\mu^2 b_4^2}{4\lambda^2}||\mathcal{H}\mathbf{1}_N\otimes v||^2
\end{split}\nonumber
\end{eqnarray}
where $0<\lambda<b_3$ is a positive constant. Under Assumption \ref{A1}, $v(t)$ is bounded for any $v(0)$ over the interval $t \in [0, \infty)$. By the comparison lemma, it can be shown that $z(t)$ is also bounded for any initial conditions over the interval $t \in [0, \infty)$.

We consider the following coordinate transformation:
\begin{eqnarray}
\begin{split}
\xi_i = z_i-v, \,\,i=1,\dots,N.
\end{split}\nonumber
\end{eqnarray}
Let $\xi = \textnormal{col}(\xi_1,\dots,\xi_N)$, we have $\xi = z - \mathbf{1}_N \otimes v$. Denote
\begin{eqnarray}
\begin{split}
\tilde{\bf{a}}(z,v)
&=(I_N \otimes M)\xi+\mathcal{T}(z)z-\mathbf{1}_N \otimes(T(v)v).
\end{split}\nonumber
\end{eqnarray}
Based on (\ref{closedz}), the closed-loop system can be written as follows:
\begin{eqnarray} \label{closedxi}
\begin{split}
\dot{\xi}=&\dot{z}-\mathbf{1}_N \otimes \dot{v}\notag\\
=&\tilde{\bf{a}}(z,v)-\mu(\mathcal{H}\otimes I_m)\xi.
\end{split}
\end{eqnarray}
Define the similar Lyapunov candidate function:
\begin{eqnarray}
W(\xi)=\xi^T(P \otimes I_m)\xi\nonumber
\end{eqnarray}
which is positive definite with $P>0$.
The derivative of $W(\xi)$ along (\ref{closedxi}) is given as follows:
\begin{eqnarray}
\begin{split}
\dot{W}(\xi) 
&=-\mu \xi^T(Q \otimes  I_m)\xi + 2\xi^T(P \otimes M)\xi\\
&+2\xi^T(P \otimes I_m)(\mathcal{T}(z)z-\mathbf{1}_N \otimes(T(v)v))\\
&\leq -\mu \lambda_{min}(Q) ||\xi||^2+2\lambda_{max}(P\otimes M) ||\xi||^2\\
&+2||\xi||||P \otimes I_m||||\mathcal{T}(z)z-\mathbf{1}_N \otimes(T(v)v)||.
\end{split}\nonumber
\end{eqnarray}
It can be verified that
\begin{eqnarray}
\begin{split}
d_j(\xi_i+v)(\xi_i+v)-d_j(v)v=d_j(z_i)z_i-d_j(v)v
\end{split}\nonumber
\end{eqnarray}
is continuously differentiable with $\zeta_j(0,v)=0$. By Lemma \ref{l3}, there is a smooth nonnegative function $\delta_j(\xi,v)$ such that
\begin{eqnarray}
\begin{split}
|d_j(\xi_i+v)(\xi_i+v)-d_j(v)v|\leq\delta_j(\xi_i,v)||\xi_i||.
\end{split}\nonumber
\end{eqnarray}
By Assumption \ref{A2},  $T(v)=\textnormal{diag}(d_1(v),\dots,d_m(v))$ is a diagonal matrix, then
\begin{eqnarray}
\begin{split}
&||\mathcal{T}(z)z-\mathbf{1}_N \otimes(T(v)v)||\\
 &\leq \sqrt{Nm}\max_{i=1,\dots,m,j=1,\dots,N}|d_j(z_i)z_i-d_j(v)v|\\
 &\leq \sqrt{Nm}\max_{i=1,\dots,m,j=1,\dots,N}\delta_j(\xi_i,v)||\xi_i||\\
 &\leq \sqrt{Nm}\beta ||\xi||
\end{split}\nonumber
\end{eqnarray}
where $\beta=\max_{i=1,\dots,m,j=1,\dots,N}\delta_j(\xi_i,v)$. This yields the following inequality
\begin{eqnarray}
\begin{split}
\dot{W}(\xi) &\leq -\mu \lambda_{min}(Q) ||\xi||^2+2\lambda_{max}(P\otimes M) ||\xi||^2\\
&+2\sqrt{Nm}\beta||P \otimes I_m||||\xi||^2.
\end{split}\nonumber
\end{eqnarray}
Next, we let
\begin{eqnarray}
\begin{split}
\mu_2 = \frac{2(\lambda_{max}(P\otimes M)+\sqrt{Nm}\beta||P \otimes I_m||)}{\lambda_{min}(Q)}.
\end{split}\nonumber
\end{eqnarray}
As a result, there exists a positive constant $\mu_2$ such that, for all $\mu > \mu_2$ and any initial conditions, the solution to the (\ref{do2}) exists and satisfies
\begin{eqnarray}
\lim_{t \rightarrow \infty} (z(t)-\mathbf{1}_N \otimes v(t)) = 0.
\end{eqnarray}
Hence,
\begin{eqnarray}
\lim_{t \rightarrow \infty} (\eta(t)-(\Phi\mathbf{1}_N) \otimes v(t)) = 0\nonumber
\end{eqnarray}
meaning that $\lim_{t \rightarrow \infty} (\eta_i(t)-\phi_i v(t))=0$, $i=1,\dots,N$.
\QEDB

\section{An Application to the Leader-Following Bipartite Consensus for Nonlinear MASs}\label{application}

In this section, we apply the NDO designed in Section \ref{ndo} to solve the BORP for a class of nonlinear MASs over static signed communication graphs.

For comparison, let us consider a class of multiple nonlinear MASs studied in \cite{liu2019distributed}:
\begin{eqnarray} \label{plant}
\begin{split}
\dot{x}_{si}&=x_{(s+1)i}, \quad s=1,\cdots,r-1\\
\dot{x}_{ri}&=f_i(x_i,v)+u_i\\
y_i &=x_{1i}\\
e_i&=y_i-\phi_iy_0, \quad i=1,\cdots,N\\
\end{split}
\end{eqnarray}
where $x_i =\textnormal{col}(x_{1i},\cdots,x_{ri})\in \mathbb{R}^{r}$ is the state for agent $i$, $u_i \in \mathbb{R}$ is the control input, $y_i \in \mathbb{R}$ denotes the  output and  $e_i \in \mathbb{R}$ is the tracking error. The functions $f_i(\cdot)$ are globally sufficiently smooth satisfying $f_i(0,0)=0$. The exogenous signals $v \in \mathbb{R}^{m}$ represent the tracking reference and/or external disturbances generated by (\ref{exo}). The functions $\phi_i$ are such that $\phi_i \in \{1,-1\}$ for $i=1,\dots,N$ and $\phi_i = 1$ if $i \in \mathcal{V}_1$, $\phi_i = -1$ if $i \in \mathcal{V}_2$.

We consider a class of distributed control law given by (\ref{u}) and formulate the leader-following bipartite consensus problem for a class of nonlinear MASs (\ref{plant}) as follows.

\begin{problem}[Leader-following Bipartite Consensus for a Class of Nonlinear MASs]
Given the nonlinear MASs (\ref{plant}) with the nonlinear exosystem (\ref{exo}) and a signed communication digraph $\bar{\mathcal{G}}^s$, design a distributed control law $u_i$ of the form (\ref{u}) such that,  for any initial states,\\
1) the trajectories of all states of the closed-loop system exist are bounded for all $t \geq 0$;\\
2) the tracking error for each agent satisfies
\begin{eqnarray}
\begin{split}
 \lim_{t \rightarrow \infty} e_i(t)=0,\,\,i=1,\dots,N.
 \end{split}\nonumber
 \end{eqnarray}
\end{problem}

 Using the NDO (\ref{do2}), we propose the following distributed control law to solve {\em Problem 2}:
\begin{eqnarray} \label{u2}
\begin{split}
u_i&=\sum_{s=1}^{r+1}\phi_i\beta_s{\bf{x}}_s(\phi_i\eta_i)-f_i(x_i,\phi_i\eta_i)-\sum_{s=1}^{r}\beta_sx_{si}\\
\dot{\eta}_i&=\phi_ia(\phi_i\eta_i)+\mu\big(\sum_{j=1}^{N}a_{ij}(\eta_j-\sgn(a_{ij})\eta_i)\\
&+a_{i0}(\phi_iv-\eta_i)\big)
\end{split}
\end{eqnarray}
where $\beta_i$, $i=1,\cdots,r$, are parameters to be assigned, $\beta_{r+1}=1$ and
\begin{eqnarray}
\begin{split}
{\bf{x}}_1(v) & = g(v)=y_0\\
{\bf{x}}_{s+1}(v) & = \frac{\partial {\bf{x}}_s(v)}{\partial v}a(v)
 \end{split}\nonumber
 \end{eqnarray}
 for $s=1,\cdots,r$. Let ${\bf{x}}(v)=\textnormal{col}({\bf{x}}_1(v),\cdots, {\bf{x}}_r(v))$.

The main result can be summarized as follows:
\begin{theorem}\label{T1}
Under Assumptions \ref{A1}-\ref{A3}, {\em Problem 2} can be solved by the distributed control law (\ref{u2}) for any $\mu>\mu_0$ where $\mu_0>0$ is some positive constant.
\end{theorem}

\emph{Proof:} Apply the distributed control law (\ref{u2}) to (\ref{plant}). This yields the closed-loop system:
\begin{eqnarray}
\begin{split}
\dot{x}_{si}&=x_{(s+1)i}, \quad s=1,\cdots,r-1\\
\dot{x}_{ri}&=f_i(x_i,v)+\sum_{s=1}^{r+1}\phi_i\beta_s{\bf{x}}_s(\phi_i\eta_i)\\
&-f_i(x_i,\phi_i\eta_i)-\sum_{s=1}^{r}\beta_sx_{si}.
\end{split}\nonumber
\end{eqnarray}
Let $x_i=\textnormal{col}(x_{1i},\cdots, x_{ri})$, and perform the coordinate transformation:
$\bar{x}_i = x_i-\phi_i{\bf{x}}(v)$
where $\bar{x}_i=\textnormal{col}(\bar{x}_{1i},\cdots, \bar{x}_{ri})$.
Then we obtain
\begin{eqnarray}
\begin{split}
\dot{\bar{x}}_{si}&=\bar{x}_{(s+1)i}, \quad s=1,\cdots,r-1\\
\dot{\bar{x}}_{ri}
&=-\sum_{s=1}^{r}\beta_s\bar{x}_{si}+f_i(x_i,v)-f_i(x_i,\phi_i\eta_i)\\
&+\sum_{s=1}^{r+1}\phi_i\beta_s{\bf{x}}_s(\phi_i\eta_i)-\sum_{s=1}^{r+1}\phi_i\beta_s{\bf{x}}_s(v)).
\end{split}\nonumber
\end{eqnarray}

The compact form can be given as follows:
\begin{eqnarray}
\dot{\bar{x}}_i = A\bar{x}_i+B\bar{u}_i, \quad i=1,\cdots,N
\end{eqnarray}
where
\begin{eqnarray}
\begin{split}
A& = \left[ \begin{array}{cc}
   0  & I_{r-1}\\
   -\beta_1& -\beta_2 \cdots -\beta_r\\
\end{array} \right], B =\textnormal{col}(0,\cdots,0,1)\\
\bar{u}_i&=f_i(x_i,v)-f_i(x_i,\phi_i\eta_i)\\
&+\sum_{s=1}^{r+1}\phi_i\beta_s{\bf{x}}_s(\phi_i\eta_i)-\sum_{s=1}^{r+1}\phi_i\beta_s{\bf{x}}_s(v))
\end{split}\nonumber
\end{eqnarray}

According to Lemma \ref{l5}, under Assumptions \ref{A1}-\ref{A3}, we attain
\begin{eqnarray}
\lim_{t \rightarrow \infty} (\eta_i(t)-\phi_i v(t)) = 0,\quad i=1,\cdots,N.\nonumber
\end{eqnarray}
exponentially for $\mu>\mu_0$ where $\mu_0$ is some positive constant. Moreover, the functions $f_i$ and ${\bf{x}}_s$ are all continuous. We conclude that $\lim_{t \rightarrow \infty} \bar{u}_i(t)=0$.

Furthermore, $A$ is a companion matrix and the parameters $\beta_i$, $i=1,\cdots,r$, are chosen such that $A$ is Hurwitz, we conclude that
\begin{eqnarray}
\lim_{t \rightarrow \infty} \bar{x}_i(t) = 0,\quad i=1,\cdots,N.
\end{eqnarray}

Hence, we can conclude that
\begin{eqnarray}
\begin{split}
\lim_{t \rightarrow \infty} e_i(t) 
&= \lim_{t \rightarrow \infty}(x_{1i}(t)-\phi_i{\bf{x}}_1(v(t))) \\
&= \lim_{t \rightarrow \infty}\bar{x}_{1i}(t)= 0,\quad i=1,\cdots,N.
\end{split}\nonumber
\end{eqnarray}
\QEDB

\section{Example}\label{exam}
\subsection{Example 1: Multiple Inverted Pendulums}
Let us consider a class of nonlinear MASs as given in \cite{slotine1991applied}:
\begin{eqnarray} \label{pendulum}
\begin{split}
\dot{x}_{1i} &= x_{2i}\\
\dot{x}_{2i} &= -2\sin(x_{1i})-x_{2i}+u_i\\
y_i&=x_{1i}\\
e_i&=y_i-\phi_iy_0,\quad i=1,2,3,4
\end{split}
\end{eqnarray}
where $x_{1i} \in \mathbb{R}$ and $x_{2i} \in \mathbb{R}$ are the state representing the position and the velocity of each inverted pendulum, $u_i \in \mathbb{R} $ is the control torque applied to each follower system, $y_i\in \mathbb{R}$ and $e_i$ are the output and the tracking error, respectively. Moreover, $x_{11}(0)=0.3, x_{21}(0)=0.4,x_{12}(0)=0.5, x_{22}(0)=0.6,x_{13}(0)=0.7, x_{23}(0)=0.8,x_{14}(0)=0.9, x_{24}(0)=1.0$.

$v=\textnormal{col}(v_1,v_2)$ can be generated by the exosystem (\ref{exo}):
\begin{eqnarray}
\begin{split}
\dot{v}_1 &= v_2\\
\dot{v}_2 &= -v_1+(1-v_1^2)v_2\\
y_0 &= v_1
\end{split}\nonumber
\end{eqnarray}
where $v_1(0)=0.1$ and $v_2(0)=0.2$. 

It can be verified that Assumptions \ref{A1} and \ref{A2} are satisfied with $M = \left[ \begin{array}{cc} 0 & 1  \\ -1 & 1  \end{array} \right]$ and $T(v)=\left[ \begin{array}{cc} 0 & 0  \\ 0 & -v_1^2  \end{array} \right]$. Moreover, the reference signal is given by $y_0=g(v)=v_1$.

\begin{figure}[H]
\centering
\resizebox*{4cm}{!}{\includegraphics{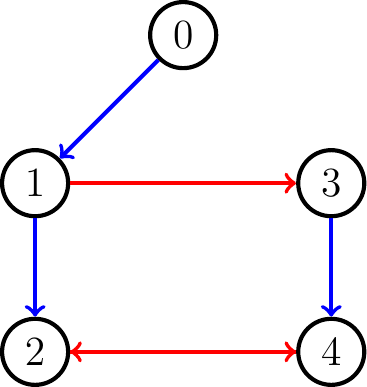}}\hspace{5pt}
\caption{Communication digraph $\bar{\mathcal{G}}^s$} \label{graph}
\end{figure}

The communication digraph is given in Fig. \ref{graph}. It is assumed that $a_{ij} \in \{-1,0,1\}$, w.l.o.g., where $a_{ij}=1$ represents the cooperation relationship denoted by the blue line while $a_{ij}=-1$ is the competition relationship denoted by the red line.

It can be seen from Fig. \ref{graph} that the communication digraph satisfies Assumption \ref{A3}. In particular, $\mathcal{V}_1 = \{1, 2\}$ and $\mathcal{V}_2 = \{3, 4\}$. Thus,
\begin{eqnarray*}
\mathcal{H}^s =\mathcal{L}^s+\Delta= \left[ \begin{array}{cccc} 1 & 0  & 0 & 0 \\ -1 &2 & 0&1\\  1& 0& 1& 0 \\0& 1& -1& 2 \end{array} \right]
\end{eqnarray*}
where the  eigenvalues of the matrix  $\mathcal{H}^s$ are $\{3,1,1,1\}$.

\begin{figure}[H]
\centering
\resizebox*{8.5cm}{!}{\includegraphics{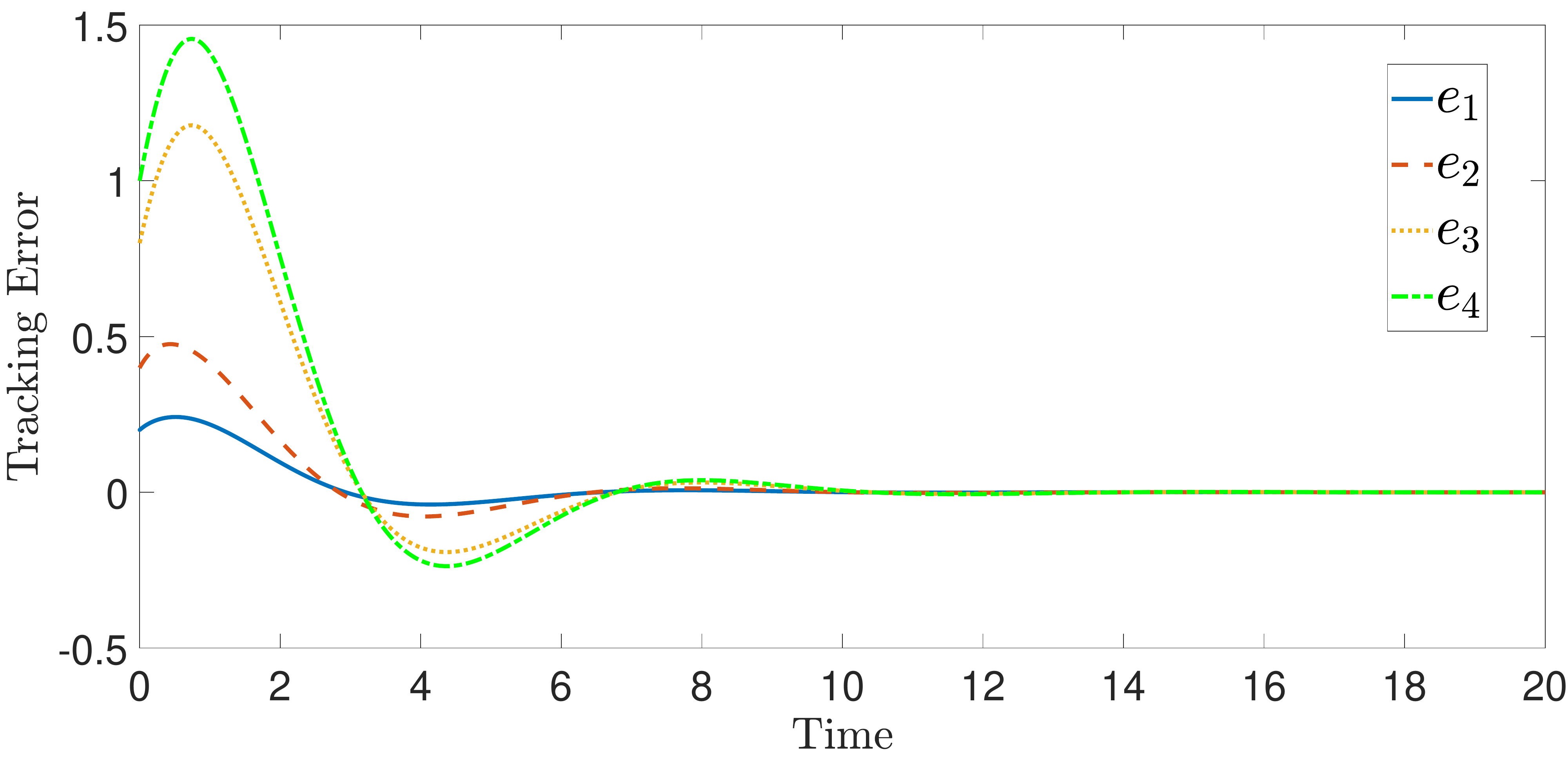}}\hspace{5pt}
\caption{Tracking error} \label{e}
\resizebox*{8.5cm}{!}{\includegraphics{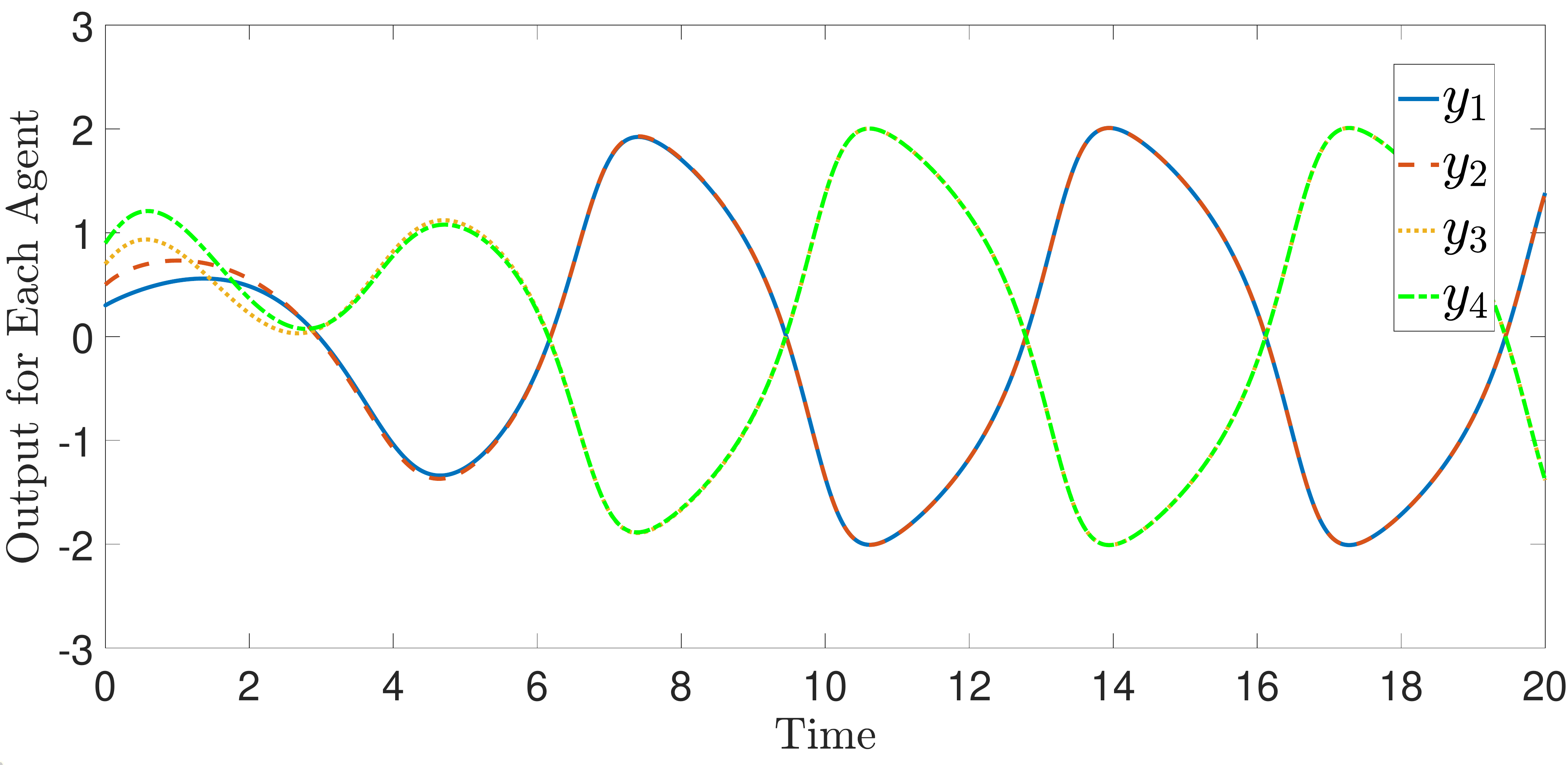}}\hspace{5pt}
\caption{Output for each agent} \label{y}
\end{figure}

Simple calculation shows that ${\bf{x}}_1(v)=v_1, {\bf{x}}_2(v) =v_2, {\bf{x}}_3(v) =-v_1+(1-v_1^2)v_2$. We apply the distributed control law (\ref{u2}) to the multiple pendulum systems (\ref{pendulum}) where the initial states of the NDO are set by $\eta_i=[(i-1)* 0.2,(i-1)*0.2]^T$. The experimental results are shown from Fig. \ref{e} to Fig. \ref{y}. It can be seen from Fig. \ref{y} that the tracking error of each agent tends to the origin as time tends to infinity, and the bipartite synchronization property can be achieved.

\subsection{Example 2: Applications to the {T}uring Pattern}

In this part, as a possible application, we apply the result in Theorem \ref{T1} to draw the zebra-stripe, a class of typical {T}uring patterns. Consider the $N=500*500 = 250000$ agents; each agent is located at one pixel of a $500*500$ figure. The dynamics of each agent are given by (\ref{pendulum}).

The signal $v= \textnormal{col}(v_1, v_2)$ can be generated by the linear exosystem $\dot{v}_1 = 0$ and $\dot{v}_2 = 0$ where $v_1(0)=-1$ and $v_2(0)=1$. In this case, the value $-1$ represents the dark colour while the value $1$ represents the white colour.

%

\begin{figure}[H]
\centering
\resizebox*{8.5cm}{!}{\includegraphics{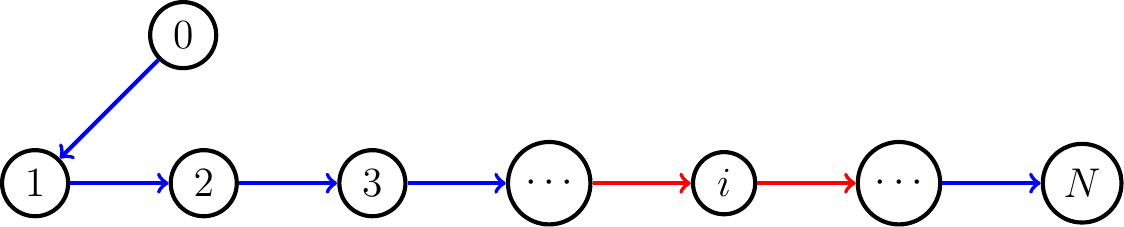}}\hspace{5pt}
\caption{The communication topology for each agent in the same line} \label{path}
\end{figure}

The communication topology for each agent in the same line is provided by Fig. \ref{path}. Node $i$ belongs to the set $V_1$ if the pixel associated with node $i$ is dark, and node $i$ belongs to the set $V_2$ if the pixel is white. If node $i$ and node $j$ belongs to the same set, then communication weight $a_{ij}$ equals to $1$ (blue color), otherwise, $a_{ij}=-1$ (red color).

\begin{figure}[H]
\centering
\resizebox*{8.5cm}{!}{\includegraphics{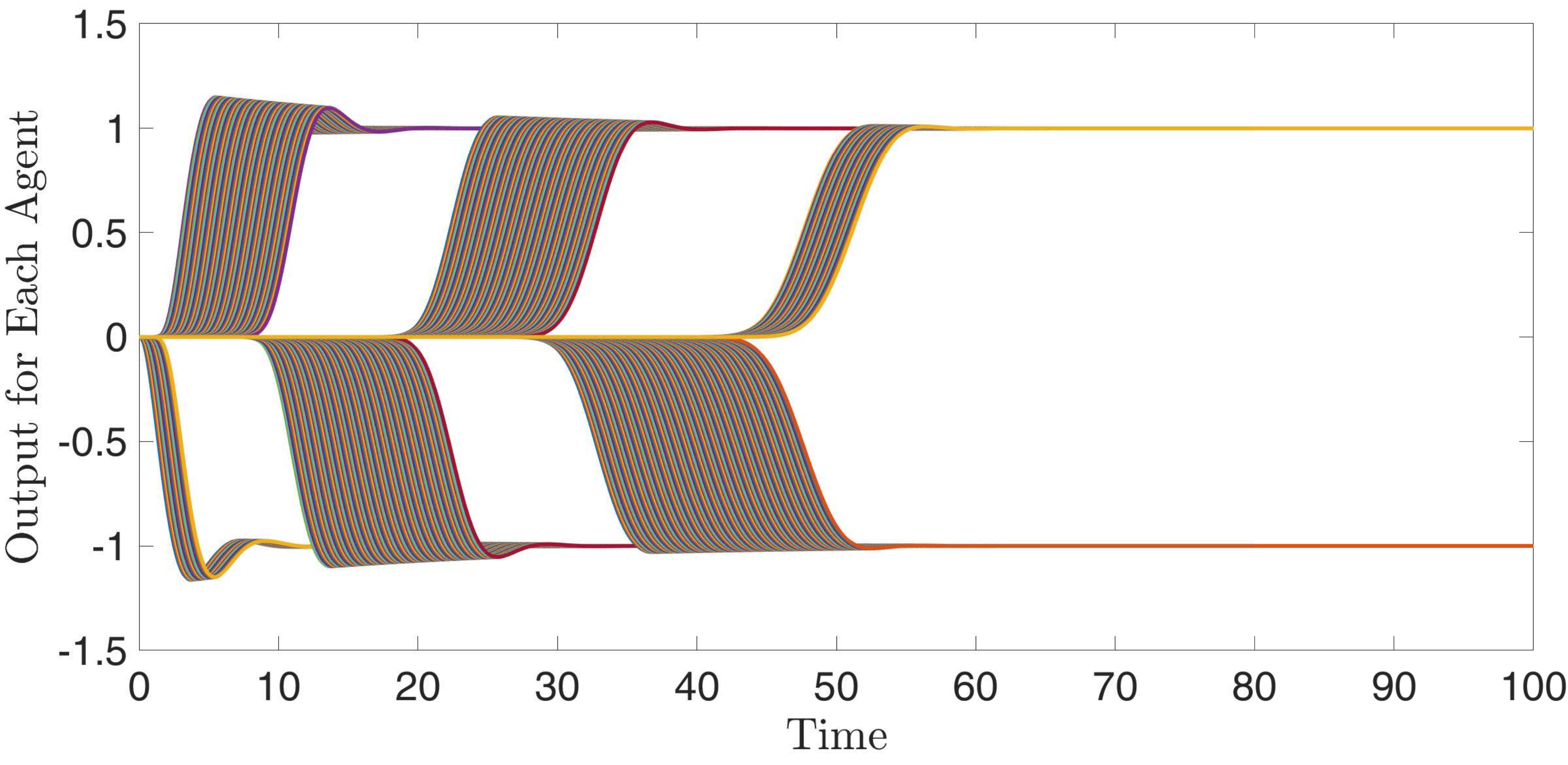}}\hspace{5pt}
\caption{Output of each agent for line 1} \label{binary}
~\\
\centering
\resizebox*{5cm}{!}{\includegraphics{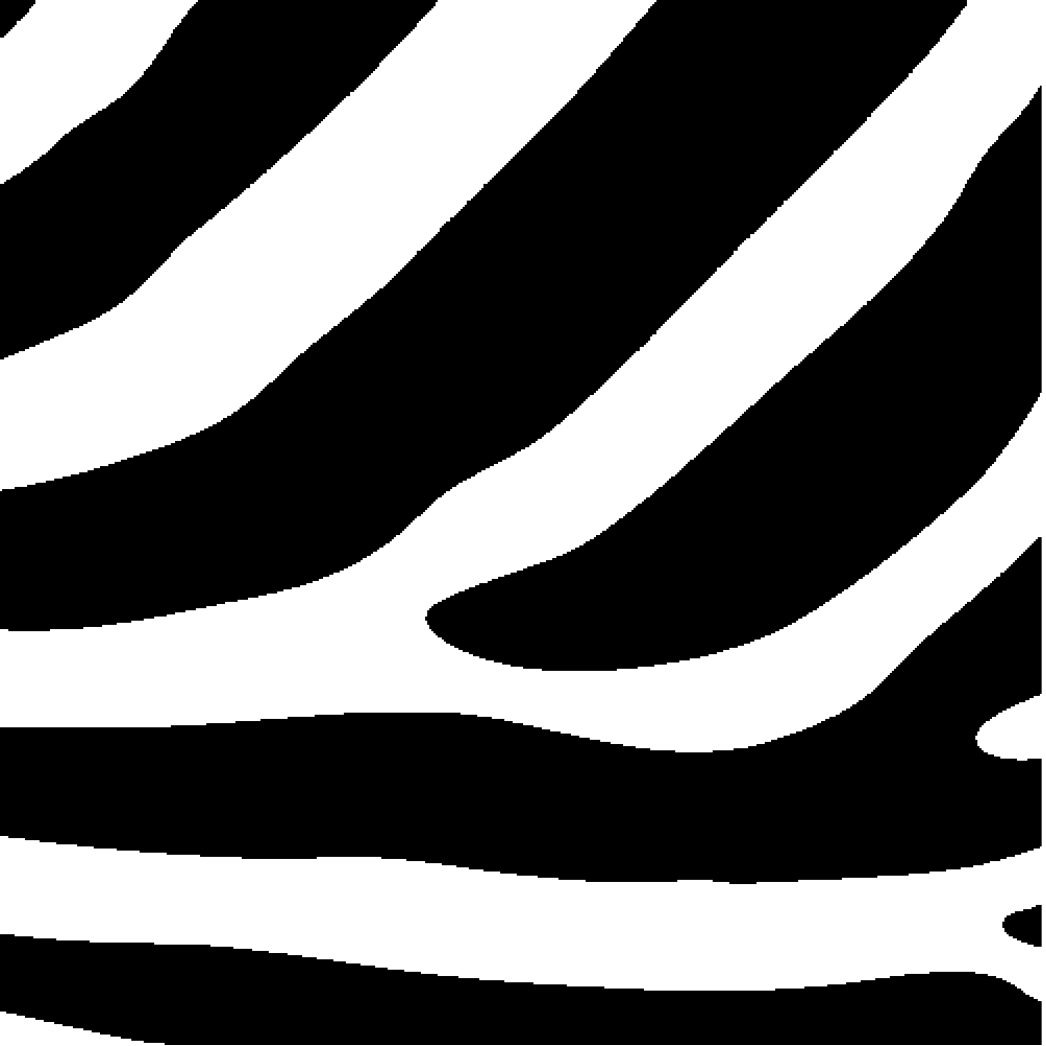}}\hspace{5pt}
\caption{Turing pattern of the zebra-stripe} \label{turing}
\end{figure}

In line 1, the output for each agent is plotted in Fig. \ref{binary}. It can be observed that some agents' outputs tend to be 1 (white colour) while the other agent's outputs tend to be -1 (black colour). Similarly, for line $i$, $i=2,\cdots,500$, each agent has its own output (-1 or 1) associated with a colour (black or white). The final result is shown in Fig. \ref{turing}, which is a typical Turing Pattern and similar to the zebra-stripe.


\begin{remark}
{It should be noted that the dynamical systems  have been used to study the {T}uring patterns with unsigned directed graphs \cite{asllani2014theory}. As the generalization of directed graphs, we have related signed graphs to the  {T}uring patterns like the zebra-stripe arising from nonlinear MASs.}
\end{remark}

\section{Conclusion}\label{con}
In this technical paper, the BORP has been formulated for nonlinear MASs.
A class of NDO over static signed communication networks has been proposed and sufficient conditions are provided to guarantee the existence and exponential stability.
In comparison to existing results presented in \cite{shtac12,aghbolagh2017bipartite,liu2019distributed}, the proposed observer applies to the leader system with general nonlinear dynamics over static signed communication networks.
As a practical application, a solution to the leader-following bipartite synchronization problem is given for a class of nonlinear MASs with cooperation-competition interactions.
It has been shown that the proposed framework can also be applied to generate {T}uring patterns.
Future work will focus on switching networks and the implementation of the proposed nonlinear observer to achieve  bipartite attitude synchronization for other classical nonlinear MASs such as Euler-Lagrange systems, multiple spacecraft systems and rigid-body systems.

\bibliographystyle{ieeetr}
\bibliography{myref}

\begin{thebibliography}{10}

\bibitem{altafini2013consensus}
C.~Altafini, ``Consensus problems on networks with antagonistic interactions,''
  {\em IEEE Transactions on Automatic Control}, vol.~58, no.~4, pp.~935--946,
  2012.

\bibitem{du2019edge}
M.~Du, B.~Ma, and D.~Meng, ``Edge convergence problems on signed networks,''
  {\em IEEE Transactions on Cybernetics}, vol.~49, no.~11, pp.~4029--4041,
  2018.

\bibitem{du2019edge2}
M.~Du, B.~Ma, and D.~Meng, ``Further results for edge convergence of directed
  signed networks,'' {\em IEEE Transactions on Cybernetics}, vol.~51, no.~12,
  pp.~5659--5670, 2019.

\bibitem{griffin2004cooperation}
A.~S. Griffin, S.~A. West, and A.~Buckling, ``Cooperation and competition in
  pathogenic bacteria,'' {\em Nature}, vol.~430, no.~7003, pp.~1024--1027,
  2004.

\bibitem{altafini2015predictable}
C.~Altafini and G.~Lini, ``Predictable dynamics of opinion forming for networks
  with antagonistic interactions,'' {\em IEEE Transactions on Automatic
  Control}, vol.~60, no.~2, pp.~342--357, 2014.

\bibitem{Meo15}
P.~De~Meo, E.~Ferrara, D.~Rosaci, and G.~M. Sarn{\'e}, ``Trust and compactness
  in social network groups,'' {\em IEEE Transactions on Cybernetics}, vol.~45,
  no.~2, pp.~205--216, 2014.

\bibitem{Xia16}
W.~Xia, M.~Cao, and K.~H. Johansson, ``Structural balance and opinion
  separation in trust--mistrust social networks,'' {\em IEEE Transactions on
  Control of Network Systems}, vol.~3, no.~1, pp.~46--56, 2015.

\bibitem{yaghmaie2017bipartite}
F.~A. Yaghmaie, R.~Su, F.~L. Lewis, and S.~Olaru, ``Bipartite and cooperative
  output synchronizations of linear heterogeneous agents: A unified
  framework,'' {\em Automatica}, vol.~80, pp.~172--176, 2017.

\bibitem{Zha17}
H.~Zhang and J.~Chen, ``Bipartite consensus of multi-agent systems over signed
  graphs: state feedback and output feedback control approaches,'' {\em
  International Journal of Robust and Nonlinear Control}, vol.~27, no.~1,
  pp.~3--14, 2017.

\bibitem{zhang2019h}
H.~Zhang, J.~Han, Y.~Wang, and H.~Jiang, ``${H}_{\infty}$ consensus for linear
  heterogeneous discrete-time multiagent systems with output feedback
  control,'' {\em IEEE Transactions on Cybernetics}, vol.~49, no.~10,
  pp.~3713--3721, 2018.

\bibitem{Jia18}
Q.~Jiao, H.~Zhang, S.~Xu, F.~L. Lewis, and L.~Xie, ``Bipartite tracking of
  homogeneous and heterogeneous linear multi-agent systems,'' {\em
  International Journal of Control}, vol.~92, no.~12, pp.~2963--2972, 2019.

\bibitem{Ma18}
C.-Q. Ma and L.~Xie, ``Necessary and sufficient conditions for leader-following
  bipartite consensus with measurement noise,'' {\em IEEE Transactions on
  Systems, Man, and Cybernetics: Systems}, vol.~50, no.~5, pp.~1976--1981,
  2018.

\bibitem{bhowmick2018leader}
S.~Bhowmick and S.~Panja, ``Leader--follower bipartite consensus of linear
  multiagent systems over a signed directed graph,'' {\em IEEE Transactions on
  Circuits and Systems II: Express Briefs}, vol.~66, no.~8, pp.~1436--1440,
  2018.

\bibitem{shi2023cooperation}
L.~Shi, Q.~Liu, J.~Shao, Y.~Cheng, and W.~X. Zheng, ``A cooperation-competition
  evolutionary dynamic model over signed networks,'' {\em IEEE Transactions on
  Automatic Control}.
\newblock doi: \url{ 10.1109/TAC.2023.3247874,2023}.

\bibitem{Liu18}
M.~Liu, X.~Wang, and Z.~Li, ``Robust bipartite consensus and tracking control
  of high-order multiagent systems with matching uncertainties and antagonistic
  interactions,'' {\em IEEE Transactions on Systems, Man, and Cybernetics:
  Systems}, vol.~50, no.~7, pp.~2541--2550, 2018.

\bibitem{shao2020}
J.~Shao, W.~X. Zheng, L.~Shi, and Y.~Cheng, ``Bipartite tracking consensus of
  generic linear agents with discrete-time dynamics over
  cooperation--competition networks,'' {\em IEEE Transactions on Cybernetics},
  vol.~51, no.~11, pp.~5225--5235, 2020.

\bibitem{jiang2017sign}
Y.~Jiang, H.~Zhang, and J.~Chen, ``Sign-consensus of linear multi-agent systems
  over signed directed graphs,'' {\em IEEE Transactions on Industrial
  Electronics}, vol.~64, no.~6, pp.~5075--5083, 2016.

\bibitem{jiang2020output}
H.~Jiang and H.~Zhang, ``Output sign-consensus of heterogeneous multiagent
  systems over fixed and switching signed graphs,'' {\em International Journal
  of Robust and Nonlinear Control}, vol.~30, no.~5, pp.~1938--1955, 2020.

\bibitem{sun2021output}
Z.~Sun, H.~Zhang, and F.~L. Lewis, ``Output sign-consensus of heterogeneous
  multiagent systems over fixed and switching signed graphs: An observer-based
  approach,'' {\em International Journal of Robust and Nonlinear Control},
  vol.~31, no.~12, pp.~5849--5864, 2021.

\bibitem{Wan17}
Y.~Wang, M.~Cheng, and H.~Zhang, ``Output bipartite consensus of heterogeneous
  linear multi-agent systems under switching topologies,'' in {\em 2017 36th
  Chinese Control Conference (CCC)}, pp.~8806--8810, IEEE, 2017.

\bibitem{aghbolagh2017bipartite}
H.~D. Aghbolagh, M.~Zamani, and Z.~Chen, ``Bipartite output regulation of
  multi-agent systems with antagonistic interactions,'' in {\em 2017 11th Asian
  control conference (ASCC)}, pp.~321--325, IEEE, 2017.

\bibitem{shtac12}
Y.~Su and J.~Huang, ``Cooperative output regulation of linear multi-agent
  systems,'' {\em IEEE Transactions on Automatic Control}, vol.~57, no.~4,
  pp.~1062--1066, 2011.

\bibitem{liang2020robust}
D.~Liang and J.~Huang, ``Robust bipartite output regulation of linear uncertain
  multi-agent systems,'' {\em International Journal of Control}, vol.~95,
  no.~1, pp.~42--49, 2022.

\bibitem{shams21hinfinity}
A.~Shams, M.~Rehan, and M.~Tufail, ``${H}_{\infty}$ bipartite consensus of
  nonlinear multi-agent systems over a directed signed graph with a leader of
  non-zero input,'' {\em International Journal of Control}, vol.~95, no.~7,
  pp.~1944--1958, 2022.

\bibitem{Qin17}
J.~Qin, W.~Fu, W.~X. Zheng, and H.~Gao, ``On the bipartite consensus for
  generic linear multiagent systems with input saturation,'' {\em IEEE
  Transactions on Cybernetics}, vol.~47, no.~8, pp.~1948--1958, 2016.

\bibitem{liang2020leader}
D.~Liang and J.~Huang, ``Leader-following bipartite consensus of multiple
  uncertain euler-lagrange systems over signed switching digraphs,'' {\em
  Neurocomputing}, vol.~405, pp.~96--102, 2020.

\bibitem{Zha16}
S.~Zhai and Q.~Li, ``Pinning bipartite synchronization for coupled nonlinear
  systems with antagonistic interactions and switching topologies,'' {\em
  Systems \& Control Letters}, vol.~94, pp.~127--132, 2016.

\bibitem{turing1990}
A.~M. Turing, ``The chemical basis of morphogenesis,'' {\em Bulletin of
  Mathematical Biology}, vol.~52, no.~1-2, pp.~153--197, 1990.

\bibitem{MaronSusan1997}
J.~L. Maron and S.~Harrison, ``Spatial pattern formation in an insect
  host-parasitoid system,'' {\em Science}, vol.~278, no.~5343, pp.~1619--1621,
  1997.

\bibitem{maini2006turing}
P.~K. Maini, R.~E. Baker, and C.-M. Chuong, ``The {T}uring model comes of
  molecular age,'' {\em Science}, vol.~314, no.~5804, pp.~1397--1398, 2006.

\bibitem{nakao2010turing}
H.~Nakao and A.~S. Mikhailov, ``Turing patterns in network-organized
  activator--inhibitor systems,'' {\em Nature Physics}, vol.~6, no.~7,
  pp.~544--550, 2010.

\bibitem{asllani2014theory}
M.~Asllani, J.~D. Challenger, F.~S. Pavone, L.~Sacconi, and D.~Fanelli, ``The
  theory of pattern formation on directed networks,'' {\em Nature
  Communications}, vol.~5, no.~1, p.~4517, 2014.

\bibitem{Su14}
Y.~Su and J.~Huang, ``Cooperative robust output regulation of a class of
  heterogeneous linear uncertain multi-agent systems,'' {\em International
  Journal of Robust and Nonlinear Control}, vol.~24, no.~17, pp.~2819--2839,
  2014.

\bibitem{liu2019distributed}
T.~Liu and J.~Huang, ``A distributed observer for a class of nonlinear systems
  and its application to a leader-following consensus problem,'' {\em IEEE
  Transactions on Automatic Control}, vol.~64, no.~3, pp.~1221--1227, 2018.

\bibitem{su2012cooperative}
Y.~Su and J.~Huang, ``Cooperative output regulation with application to
  multi-agent consensus under switching network,'' {\em IEEE Transactions on
  Systems, Man, and Cybernetics, Part B (Cybernetics)}, vol.~42, no.~3,
  pp.~864--875, 2012.

\bibitem{huang2017cooperative}
J.~Huang, ``The cooperative output regulation problem of discrete-time linear
  multi-agent systems by the adaptive distributed observer,'' {\em IEEE
  Transactions on Automatic Control}, vol.~62, no.~4, pp.~1979--1984, 2016.

\bibitem{wang2018adaptive}
S.~Wang and J.~Huang, ``Adaptive leader-following consensus for multiple
  {E}uler--{L}agrange systems with an uncertain leader system,'' {\em IEEE
  Transactions on Neural Networks and Learning Systems}, vol.~30, no.~7,
  pp.~2188--2196, 2018.

\bibitem{wang2021event}
S.~Wang, Z.~Shu, and T.~Chen, ``Event-triggered attitude synchronization of
  multiple rigid-body systems,'' {\em Systems \& Control Letters}, vol.~149,
  p.~104879, 2021.

\bibitem{michel1977qualitative}
A.~N. Michel and R.~K. Miller, ``Qualitative analysis of large scale dynamical
  systems,'' 1977.

\bibitem{chen2015stabilization}
Z.~Chen and J.~Huang, ``Stabilization and regulation of nonlinear systems,''
  {\em Cham, Switzerland: Springer}, 2015.

\bibitem{slotine1991applied}
J.-J.~E. Slotine, W.~Li, {\em et~al.}, {\em Applied nonlinear control},
  vol.~199.
\newblock Prentice Hall Englewood Cliffs, NJ, 1991.

\end{thebibliography}

\end{document}